\documentclass[11pt]{amsart}
\usepackage{amsthm}
\usepackage{amsmath}
\begin{document}
\newcommand{\Q}{{\mathbb Q}}
\newcommand{\R}{{\mathbb R}}
\newcommand{\Z}{{\mathbb Z}}
\newcommand{\F}{{\mathbb F}}
\renewcommand{\P}{{\mathbb P}}
\renewcommand{\O}{{\mathcal O}}
\newcommand{\Pic}{{\rm Pic\,}}
\newcommand{\Ext}{{\rm Ext}\,}
\newcommand{\rank}{{\rm rk}\,}
\newcommand{\sbull}{{\scriptstyle{\bullet}}}
\newcommand{\bX}{X_{\overline{k}}}
\newcommand{\ch}{\operatorname{CH}}
\newcommand{\tors}{\text{tors}}
\newcommand{\cris}{\text{cris}}
\newcommand{\alg}{\text{alg}}
\let\isom=\simeq
\let\rk=\rank
\let\tensor=\otimes

\newtheorem{theorem}{Theorem}[section]      
\newtheorem{lemma}[theorem]{Lemma}          %
\newtheorem{corollary}[theorem]{Corollary}  
\newtheorem{proposition}[theorem]{Proposition}

\theoremstyle{definition}
\newtheorem{conj}[theorem]{Conjecture}
\newtheorem*{example}{Example}
\newtheorem{question}[theorem]{Question}

\theoremstyle{definition}
\newtheorem{remark}[theorem]{Remark}

\numberwithin{equation}{section}

\title{Exotic Torsion, Frobenius splitting and the slope spectral sequence}
\author{Kirti Joshi}
\address{Math. department, University of Arizona, 617 N Santa Rita, Tucson
85721-0089, USA.}

\maketitle

\tableofcontents

\section{Introduction}
   In this paper we continue our investigations on some questions of
V.~B.~Mehta on the crystalline aspects of Frobenius splitting which
were begun in \cite{joshi98}.  In response to a question of Mehta, we
show that any $F$-split smooth projective threefold is Hodge-Witt (see
Theorem~\ref{frob-hodge-witt}).  After the first draft of this paper
and \cite{joshi98} were written, Rajan and the author showed that
there exist examples of Frobenius split varieties which are not
ordinary in the sense of Bloch-Kato-Illusie-Raynaud (see
\cite{bloch86}, \cite{illusie83b}) and if the dimension is bigger than
four then these examples are not even Hodge-Witt (these examples are
now included in \cite{joshi98}). Thus the result of
Theorem~\ref{frob-hodge-witt} on the Hodge-Witt property of Frobenius
split threefolds is best possible. Frobenius split surfaces are
ordinary (see \cite{joshi98}) and hence Hodge-Witt.  The method of
proof uses \cite{illusie83b} and \cite{ekedahl84} to extract an
explicit criterion for the degeneration of the slope spectral sequence
of any smooth projective threefold (see Theorem~\ref{criterion}). This
criterion can be applied to Frobenius split varieties via
Theorem~\ref{finiteness} (see \cite{joshi98}). This criterion also
applies to any unirational or Fano threefolds and we deduce that these
are Hodge-Witt (see Theorem~\ref{unirat-fano}).

Mehta has also raised the following question.
\begin{question}\label{mehtaq}
   Suppose $X$ is a smooth, projective and $F$-split variety,
then does there exist a finite \'etale cover $X'\to X$ such that
$\Pic(X')$ is reduced and $H^2_{\rm cris}(X'/W)$ is torsion free?
\end{question}

   In an attempt to answer the above question of Mehta we show that
the second crystalline cohomology of a Frobenius split variety has no
exotic torsion (see Theorem~\ref{exotic}). This result is the best
possible as a variant of an example of Igusa (see \cite{igusa55})
shows. Igusa's construction
gives an example of an $F$-split variety whose Picard scheme is not
reduced; indicating the presence of torsion in the second crystalline
cohomology. Our method of proof also shows, more generally, that the
vanishing of $H^1(X,B_1\Omega^1_X)$ is sufficient to ensure the
vanishing of exotic torsion in the second crystalline cohomology of
$X$. 

      Theorem~\ref{criterion} has other applications when combined
with \cite{suwa88}, \cite{gros88b} and \cite{joshi98}. We deduce that
the $p$-torsion in the Chow group of codimension two cycles on any
smooth projective Frobenius split threefold is of finite cotype (i.e.,
the Pontryagin dual of the $p$-torsion is a direct sum of a finite
number of copies of $\Z_p$ and a finite group). As another application
of Theorem~\ref{finiteness}, this time to $\ell$-adic Abel-Jacobi
mappings of \cite{bloch79}, we note that the restriction of the
Abel-Jacobi map (see \cite{bloch79}) to the group of $\ell$-primary
torsion codimension $i$ cycles (for $i\geq 1$) which are algebraically
equivalent to zero (modulo rational equivalence) on any projective,
smooth Frobenius split variety is not surjective if
$H^{2i-1}(X,W(\O_X))\tensor W[1/p]$ is non-zero.  This applies in
particular to odd dimensional Frobenius split Calabi-Yau varieties
(see Corollary~\ref{calabi-yau}).

   We also include here an unpublished result of Mehta on
liftability to $W_2$ of Frobenius split varieties and its application
via \cite{deligne87} to the degeneration of the Hodge to de Rham.

   Our debt of gratitude to Minhyong Kim, V.~B.~Mehta and
C.~S.~Rajan is perhaps too extensive to document; their comments on
numerous proto-versions of this paper have gone a long way in shaping
the present version. Thanks are also due to Douglas Ulmer for his
constant support and encouragement during the course of this work. We
are also grateful to Luc Illusie for his comments.

      Throughout this paper following notations will be in force. Let
$k$ be an algebraically closed field of characteristic $p>0$, $W=W(k)$
is ring of Witt vectors of $k$, for all $n\geq 1$, we write
$W_n(k)=W/p^n$.  Let $K$ be the quotient field of $W(k)$.

\section{Frobenius split varieties}
   Let $X/k$ be a smooth, projective variety over an
algebraically closed field $k$ of characteristic $p>0$.  Let $F:X\to
X$ be the absolute Frobenius morphism of $X$. Following Mehta and
Ramanathan (see \cite{mehta85}), we say that $X$ is Frobenius split
(or simply $F$-split) if the exact sequence
\begin{equation}
0 \to \O_X \to F_*(\O_X)
\to B_1\Omega^1_X \to 0,
\end{equation}
 splits as a sequence of locally
free $\O_X$ modules.  

   This notion was defined in \cite{mehta85}, where a number of
principal properties of such varieties were proved.  The notion of
Frobenius splitting has played an important role in the study of
Schubert varieties and algebraic groups since its inception in
\cite{mehta85}.  In \cite{mehta87}, it was shown that an Abelian
variety is $F$-split if and only if it is ordinary. In \cite{joshi98},
it was shown that any smooth, projective $F$-split surface is
Block-Kato ordinary (for this notion see \cite{bloch86},
\cite{illusie83b}).
  
   We will say that a smooth projective variety $X$ is a
Calabi-Yau variety if the canonical bundle of $X$ is trivial and
$H^i(X,\O_X)=0$ for $0<i<\dim(X)$. From \cite{mehta85} it follows that
a Calabi-Yau variety of dimension $n=\dim(X)$ is Frobenius split if
and only if $F:H^n(X,\O_X)\to H^n(X,\O_X)$ is injective.

\section{The slope spectral sequence}
   We recall the formalism of de Rham-Witt complexes as developed
in \cite{bloch86}, \cite{illusie79b} and
\cite{illusie83b}. Specifically, let $R$ be the
Cartier-Dieudonne-Raynaud ring over $k$, recall that this is a graded
$W$-algebra generated by symbols $F,V,d$, which is generated in degree
zero by $F,V$ with relations $FV=VF=p$ and $Fa=\sigma(a)F$ for all
$a\in W$, and $aV=V\sigma^{-1}(a)$ for all $a\in W$; and $R$ is
generated in degree $1$ by $d$ with properties $FdV=d$ and $d^2=0$,
$da=ad$ for all $a\in W$. Thus $R=R^0\oplus R^1$ where $R^0$ is the
usual Cartier-Dieudonne ring and $R^1$ is an $R^0$-bimodule generated
by $d$.  A graded module over $R$ is a complex of $R^0$-modules
$\cdots\to M^i \to M^{i+1}\to\cdots$ where the differentials satisfy
$FdV=d$.

   The de Rham-Witt complex on $X$ is a sheaf of graded
$R$-modules on $X$, and is denoted by $W\Omega^\sbull_X$, and one has
$W\Omega^0_X=W(\O_X)$ is the sheaf of Witt vectors on $X$ constructed
by Serre (see \cite{serre58}). The cohomology of the de Rham-Witt
complex, denoted $H^j(X,W\Omega^\sbull_X)$ in the sequel, inherits the
structure of a graded $R$-module. As was shown in \cite{illusie79b},
\cite{bloch86}, the de Rham-Witt complex computes the crystalline
cohomology of $X$ and we also have a spectral sequence (the slope
spectral sequence):
\begin{equation}
E_1^{i,j}=H^j(X,W\Omega^i_X) \implies H^{i+j}_{\rm cris}(X/W)
\end{equation}
The groups $H^j(X,W\Omega^i_X)$ are not always of finite type over $W$
and we say that $X$ is Hodge-Witt if for all $i,j$, the groups
$H^j(X,W\Omega^i_X)$ are finite type modules over $W$.

\section{Dominoes and Ekedahl's duality}
   We recall briefly the formalism of dominoes developed in
\cite{illusie83b}, \cite{ekedahl84}. Let $M$ be a graded $R$-module
over the Cartier-Dieudonne-Raynaud ring; fix and $i\in \Z$ and we let
$V^{-r}Z^i=\{x\in M^i| V^r(x)\in Z^i\}$ for any $ r\geq 0$, and where
$Z^i=\ker(d:M^i\to M^{i+1})$. These form an decreasing sequence of
$W$-submodules of $M^i$ and we let $V^{-\infty}Z^i=\cap_{r\geq 0}
V^{-r}Z^i$. Similarly let $F^\infty B^{i+1}=\cup_{s\geq 0}
F^s(B^{i+1})$ where $B^{i+1}=image(d:M^i\to M^{i+1})$. The
differential $d:M^i\to M^{i+1}$ factors canonically as
\begin{equation}
M^i\to M^i/V^{-\infty}Z^i \to F^\infty B^{i+1} \to M^{i+1}.
\end{equation}

   Let $M=M^0\to M^1$ be a graded $R$-module concentrated in two
degrees. We say $M$ is a dominoe if $V^{-\infty}Z^0=M^0$ and $F^\infty
B^1=M^1$. Any dominoe has a filtration by sub-$R$-modules (graded)
such that graded pieces are certain standard dominoes (see
\cite{illusie83b}). The length of the filtration is an invariant of
the dominoe and is equal to $\dim_kM^0/VM^0$ and is called the
dimension of the dominoe $M$.

   Let $M$ be a graded $R$-module. For all $i$, the standard
factorization of the differential $M^i\to M^i/V^{-\infty}Z^i \to
F^\infty B^{i+1} \to M^{i+1}$ gives an associated dominoe
$M^i/V^{-\infty}Z^i \to F^{\infty} B^{i+1}$.  By Proposition~2.18 of
\cite{illusie83b} we have: the differential $d:M^i\to M^{i+1}$ is zero
if and only if the associated dominoe is zero dimensional.

   In particular, we will write $T^{i,j}$ for the dimension of
the dominoe associated to the differential $d:H^j(X,W\Omega_X^i) \to
H^j(X,W\Omega^{i+1}_X)$.

   The key input for our result is the following duality theorem
due to T.~Ekedahl (see \cite{ekedahl84}).

\begin{theorem}[Ekedahl's duality for dominoes]\label{ekedahlduality}
   Let $X$ be a smooth, projective variety of dimension $n$ over
a perfect field $k$ and let $T^{i,j}$ be the dimension of the dominoe
associated to the differential $H^j(W\Omega^i_X) \to
H^j(W\Omega_X^{i+1})$. Then we have
\begin{equation}
\label{duality}
   T^{i,j}=T^{n-i-2,n-j+2}
\end{equation}
\end{theorem}

\section{A finiteness result}
   Let $X/k$ be a smooth projective, Frobenius split variety over
$k$. The following theorem was proved in \cite{joshi98}. We give a
proof here, of this fact, from a slightly 
different point of view. We will need this result throughout 
this paper.

\begin{theorem}\label{finiteness}
	Let $X$ be as above. Then for all $i\geq 0$ the Hodge-Witt
groups $H^i(X,W(\O_X))$ are of finite type $W$-modules.
\end{theorem}
\begin{proof}
We first show that the differential $H^i(X,W(\O_X))\to
H^i(X,W\Omega^1_X)$ is zero. To prove this, we observe that by
Corollary~3.8 of \cite{illusie83b} it suffices to show that
$H^i(X,W(\O_X))$ is $F$-finite. Observe that to prove $F$-finiteness
it suffices to prove that $H^i(X,B_n\Omega^1_X)$ is of bounded
dimension as $n$ varies. This is where we use the fact that $X$ is
$F$-split.  It is immediate (see \cite{joshi98}) from the definition
of $F$-splitting that $H^i(X,B_1\Omega^1_X) =0$ for all $i\geq 0$.  By
\cite{illusie79b} we have an exact sequence
\begin{equation}\label{devisage} 0 \to B_1\Omega^1_X \to
B_{n+1}\Omega^1_X \to B_n\Omega^1_X \to 0
\end{equation}
 so the required boundedness follows from this by induction on
$n$.

   As there are no de Rham Witt differentials which map to
$H^i(X,W(\O_X))$ we see that $E_2^{1,i}=H^i(X,W(\O_X))$. But by
\cite{illusie83b} we know that the $E_2$ terms of the slope spectral
sequence are all of finite type $W$-modules. Hence we see that
$H^i(X,W(\O_X))$ is of finite type.
\end{proof}

\section{On the slope spectral sequence of threefolds}
\label{degeneration}
   Throughout this section $X/k$ is a smooth projective, variety of
dimension three. The main aim of this section is to give a fairly
explicit criterion, which is the analogue for threefolds of
\cite{nygaard79b}, for the degeneration of the slope spectral sequence
at $E_1$.

\begin{theorem}\label{criterion}
   Let $X/k$ be a smooth, projective threefold. Then the slope
spectral sequence of $X$ degenerates at $E_1$ if and only if for all
$i\geq 0$, the groups $H^i(X,W(\O_X))$ are finite type modules over
$W$. In other words, $X$ is Hodge-Witt if and only if $H^i(X,W(\O_X))$
are of finite type for $i\geq 0$.
\end{theorem}

\begin{proof}
   By definition we know that when $X$ is Hodge-Witt then the
relevant cohomology groups are of finite type over $W$. So it suffices
to prove that $X$ is Hodge-Witt under the assumption that the
cohomology groups $H^i(X,W(\O_X))$ are of finite type over $W$.

   We note a number of reductions. To prove that these cohomology
groups are of finite type over $W$ it suffices to show that all the
differentials in the slope spectral sequence are zero: because all the
terms of the slope spectral sequence from $E_2$ onwards are all of
finite type over $W$.  Next we note that by Corollary~2.17, Page~614
of \cite{illusie79b}, $H^0(X,W\Omega^i_X)$ are finite type for all
$i\geq 0$. By Corollary~2.18, Page 614 of \cite{illusie79b} we know
that $H^j(X,W\Omega^n_X)$ are finite type for all $j\geq 0$.  By
Corollary~3.11 of \cite{illusie83b} we know that $H^1(X,W\Omega^i_X)$
are of finite type over $W$.  Thus we are left with  the following
cohomologies and differentials:
\begin{equation}
\begin{matrix}
H^2(X,W(\O_X)) & \to & H^2(X, W\Omega^1_X) & \to & H^2(X, W\Omega_X^2)\\
H^3(X,W(\O_X)) & \to & H^3(X, W\Omega^1_X) & \to & H^3(X, W\Omega_X^2)
\end{matrix}
\end{equation}

   Now by \ref{finiteness} the differentials on the left column
are zero. By Ekedahl's duality (see Theorem~\ref{ekedahlduality})
$T^{1,3}=T^{3-1-2,3-3+2}=T^{0,2}$ and
$T^{1,2}=T^{3-1-2,3-2+2}=T^{0,3}$ and as the left column of
differentials is zero, we have $T^{0,3}=T^{0,2}=0$. This completes the
proof.
\end{proof}

\begin{corollary}\label{frob-hodge-witt}
Let $X$ be a smooth, projective and $F$-split threefold. Then $X$ is
Hodge-Witt.
\end{corollary}

\begin{proof}
This is immediate from Theorem~\ref{finiteness} 
and Theorem~\ref{criterion}.
\end{proof}

The structure of crystalline cohomology of Hodge-Witt varieties has
been studied in detail in \cite{illusie83b}, and we recall the
following here for the readers convenience.

\begin{theorem}\label{newton-hodge}
Let $X$ be a smooth projective  threefold. Then the
crystalline cohomology of $X$ admits a Newton-Hodge decomposition:
\begin{equation}
H^n_{cris}(X/W)=\oplus_{i+j=n} H^j(X,W\Omega^i_X)
\end{equation}
and one has further decomposition into pure and fractional slope
parts: 
\begin{equation}
H^j(X,W\Omega^i_X)=H^{i+j}_{[i]}\oplus H^{i+j}_{]i,i+1[},
\end{equation}
where the first term has pure slope $i$ and the other term consists of
slopes between strictly between $i,i+1$.
\end{theorem}

\section{Exotic Torsion}\label{exotictorsion}
	In this section we prove that any smooth, projective,
Frobenius split variety does not have exotic torsion in its second
crystalline cohomology (for the definition of exotic torsion see
\cite{illusie79b}).

   In \cite{illusie79b} Illusie gave the following devisage of
torsion in $H^2_{cris}(X/W)$ for any smooth projective variety
$X/k$. There is a $W$-submodule of torsion in $H^2_{cris}(X/W)$ which
is called the divisorial torsion and the quotient of torsion in
$H^2_{cris}(X/W)$ by the divisorial torsion is the exotic
torsion. Examples of exotic torsion, while hard to construct, do
exist. In some sense exotic torsion is geometrically least understood
part of torsion in second crystalline cohomology. Divisorial torsion,
on the other hand has geometric manifestation: for instance, part of
divisorial torsion (the $V$-torsion) is zero if and only if the Picard
scheme is reduced, while the remaining part of divisorial torsion is
the torsion arising from the Neron-Severi group of $X$ via the
crystalline cycle class map (see \cite{illusie79b}).

	Our result on exotic torsion together with the fact that an
\'etale cover of a Frobenius split variety is again Frobenius split
shows that we may ignore exotic torsion in Mehta's
Question~\ref{mehtaq}.  We note that in \cite{mehta87} it was shown
that any smooth, projective, $F$-split variety with trivial tangent
bundle has a finite Galois \'etale cover which is an Abelian variety
and it is standard (see for instance \cite{illusie79b}) that the
crystalline cohomology of an Abelian variety is torsion free.

The main results of this section are Theorem~\ref{exotic} and
Theorem~\ref{exotic-criterion}.

\begin{theorem}\label{exotic}
Let $X/k$ be a smooth, projective, $F$-split variety over an algebraically
closed field of characteristic $p>0$. Then the second crystalline
cohomology of $X$ has no exotic torsion.
\end{theorem}
\begin{proof}
 We first show that for all
$i\geq 0$ we have
    $$F:H^i(X,W(\O_X)) \to H^i(X,W(\O_X)),$$ is an isomorphism. To see
this we proceed as follows. We have an exact sequence (see
\cite{illusie79b})
\begin{multline}    
    \cdots\to H^{i-1}(X,W(\O_X)/FW(\O_X)) \to H^i(X,W(\O_X)) \to \\
      \qquad \xrightarrow{F} H^i(X,W(\O_X)) \to 
      H^i(X,W(\O_X)/FW(\O_X))\to\cdots
\end{multline}
    So it suffices to show that $H^i(X,W(\O_X)/FW(\O_X))$ is zero. But
$$H^i(X,W(O_X)/FW(\O_X))=\varprojlim_{n} H^i(X,B_n\Omega^1_X)$$ (by
\cite{illusie79b}, page 609, 2.2.2)
But as $X$ is $F$-split by the proof Theorem~\ref{finiteness}, we know
that the groups on the right are all zero. 

    We are now ready to prove Theorem~\ref{exotic}. Recall that by
    \cite{illusie79b}
6.7.3, Page~643, exotic torsion in $H^2_{\rm cris}(X/W)$,
denoted $H^2_{\rm cris}(X/W)_e$, is the quotient 
\begin{equation}
   H^2_{\rm cris}(X/W)_e=\frac{Q^2}{H^2(X,W(\O_X))_{\rm V-tors}}
\end{equation} 
where 
\begin{equation} 
   Q^2={\rm Image}(H^2_{\rm cris}(X/W)_{\rm Tor}\to H^2(X,W(\O_X)).  
\end{equation}

   Thus to prove our result it will suffice to show that all the
$p$-torsion in $H^2(X,W(\O_X))$ is also $V$-torsion.  But we have seen
that $F$ is an automorphism of $H^2(X,W(\O_X))$ and the relation
$FV=p$ then shows that $V=F^{-1}p$.  So all the $p$-torsion is
$V$-torsion. This finishes the proof.
\end{proof}

\begin{corollary}
Let $X$ be a smooth, projective Frobenius split variety over a perfect
field $k$. Then $\Pic(X)$ is reduced if and only if
$H^2(X,W(\O_X))_{\rm Tor}= 0$ 
\end{corollary}
\begin{proof}
It is clear from \cite{illusie79b} that $\Pic(X)$ is reduced if and
only if $H^2(X,W(\O_X))$ does not contains $V$-torsion. But by the
proof of Theorem~\ref{exotic} we see that all the $p$-torsion of
$H^2(X,W(\O_X))$ is $V$-torsion, hence the assertion.
\end{proof}

The method of proof of Theorem~\ref{exotic} can be distilled to
provide the following sufficient condition for the vanishing of exotic
torsion in any smooth projective variety:

\begin{theorem}\label{exotic-criterion}
Let $X$ be a smooth, projective variety over a perfect field $k$ of
characteristic $p>0$. Suppose that $H^1(X,B_1\Omega^1_X)=0$, then
$H_{\rm cris}^2(X/W)$ does not have exotic torsion.
\end{theorem}

\begin{proof}
From the proof of Theorem~\ref{exotic}, it is clear that if $F$
is injective on $H^2(X,W(\O_X))$ then $X$ has no exotic torsion.
So it suffices to verify that the $F$-torsion is zero under our
hypothesis. By Corollary~3.5, page~569 and Corollary~3.19 of
\cite{illusie79b}, we get 
\begin{equation} 0 \to H^1(W(\O_X))/F
\to H^1(W(\O_X)/F) \to {}_FH^2(W(\O_X))\to 0 
\end{equation} 
where the $F$-torsion of $H^2(X,W(\O_X))$ (the term on the
extreme right in the above equation) is the image of
$H^1(X,W(\O_X)/FW(\O_X))$. But by Corollary~2.2, page~609 of
\cite{illusie79b}, we know that \begin{equation}
H^1(X,W(\O_X)/FW(\O_X))=\varprojlim_{n}H^1(X,B_n\Omega^1_X)
\end{equation} By our hypothesis  $H^1(X,B_1\Omega^1_X)=0$,  by
induction on $n$ using \eqref{devisage}, one sees that
$H^1(X,B_n\Omega^1)=0$ for all $n$ so that the inverse limit is
zero and hence the $F$-torsion of $H^2(X,W(\O_X))$ is zero.
\end{proof}

\begin{remark}
	It is easy to see that the exotic torsion may be zero while
$$H^1(X,B_1\Omega^1_X)\neq 0.$$ The standard example of this
phenomenon is when $X$ is a supersingular K3 surface, where one knows
that $H^1(X,B_n\Omega_X^1)$ is $n$-dimensional for $n\geq 1$, but
$H^2_{\rm cris}(X/W)$ is torsion free while $H^2(X,W(\O_X))$ is
$F$-torsion.
\end{remark}

\section{Ordinarity of the Albanese Scheme}  
   In this section we prove that the reduced Picard scheme of a
Frobenius split variety is ordinary. This result can be proved by
other methods as well, but we present a proof here using the methods 
of the previous sections. As the Albanese variety is the dual of the
reduced Picard scheme of $X$ we state the results in terms of
$\Pic^0(X)_{\rm red}$.

\begin{theorem}\label{slopezero}
   Let $X/k$ be a smooth, projective, $F$-split variety and assume
that $k$ is algebraically closed. Then for all $i\geq 0$, the natural
map
\begin{equation}
   H^i_{et}(X,\Z_p)\tensor_{\Z_p} W(k)\to H^i(X,W(\O_X))
\end{equation}
is an isomorphism.
\end{theorem}
\begin{proof}
   By Corollary~3.5 of \cite{illusie83b} we have the following
exact sequence
   \begin{equation}
    0 \to H^i(X,W\Omega^0_{\log}) \to H^i(X,W(\O_X)) \xrightarrow{1-F} 
    H^i(X,W(\O_X)) \to 0,
   \end{equation}
  and by definition (see \cite{illusie83b}, Page~191) we know that 
      \begin{equation}
        H^i(X,W\Omega^0_{\log})\isom H^i_{et}(X,\Z_p).
      \end{equation}

   Next we know by Theorem~\ref{finiteness} (see\cite{joshi98}), that
$H^i(X,W(\O_X))$ are all finite type $W$-modules.  Now from
Lemma~6.8.4, page~643 of \cite{illusie79b} and the fact that $F$ is an
automorphism on $H^i(X,W(\O_X))$ we see that the natural map
\begin{equation} \ker(1-F)\tensor_{\Z_p}
W(k)=H^i(X,\Z_p)\tensor_{\Z_p}W(k)\to  H^i(X,W(\O_X)) \end{equation}  is
an isomorphism and this finishes the proof. 
\end{proof}

\begin{theorem}\label{ordinarity}
   Let $X$ be any smooth, projective, $F$-split variety,  
let $A=\Pic(X)^0_{\rm red}$. Then $A$ is an
 ordinary Abelian variety.
 \end{theorem}

 \begin{proof} If $A=0$ then there is nothing to prove so we will
assume that $A$ is non-zero. We show that the $F$-crystal $(H^1_{\rm
cris}(X/W)\tensor K,F)$ which is the Dieudonne module of the reduced
Picard variety of $X$, has only two slopes $0,1$. This together with
duality shows that the Abelian variety $A$ is ordinary. So it remains
to check the statement about slopes.  From Theorem~\ref{slopezero} we
know that
\begin{equation}     
\ker(1-F)\tensor_{\Z_p} W(k)=H^1(X,\Z_p)\tensor_{\Z_p}W(k)\to
H^1(X,W(\O_X)) 
\end{equation} 
   is an isomorphism. As the space on the right is the part of
$H^1_{\rm cris}(X/W)$ on which Frobenius acts through slopes $[0,1[$
and the cohomology on the left is the unit root sub crystal (i.e., part
of the cohomology on which Frobenius operates via slopes zero) of
$H^1_{\rm cris}(X/W)$ (see \cite{illusie79b}, Page~627, Paragraph
5.4). Thus we see that the slope zero part of the Frobenius is
isomorphic to the part with slopes between $[0,1[$.  Thus the only
possible slopes for the $F$-crystal $(H^1_{\rm cris}(X/W),F)$ are zero
and one. This proves the assertion. 
\end{proof}

\section{Liftings to $W_2$; Hodge-de Rham spectral sequence}
    We collect together a few facts about $F$-split varieties.  The
following consequence of $F$-splitting and Deligne-Illusie  technique (see \cite{deligne87}) 
is due to V.~B.~Mehta (unpublished).

\begin{theorem}\label{mehta-degen}
   Let $X$ be a smooth, projective, $F$-split variety over an algebraically
   closed field $k$ of characteristic $p>0$.
   Then for all  $i+j<p$,
the Hodge to de Rham spectral sequence degenerates at $E_1^{i,j}$. In
particular, for $i+j=1$ we have the following exact sequence
\begin{equation}    
    0 \to H^0(X,\Omega_X^1) \to H^1_{\rm DR}(X/k) \to H^1(X,\O_X) \to
         0.
\end{equation}
    Moreover, any $F$-split variety with $\dim(X)<p$ satisfies
Kodaira-Akizuki-Nakano vanishing. \end{theorem}

\begin{proof} 
    After \cite{deligne87} we only need to check that $X$ admits
a proper flat lifting to $W_2(k)$. As was noted in \cite{srinivas90} the
obstruction to lifting a variety to $W_2(k)$ is the image of the
obstruction to lifting the pair $(X,F)$ to $W_2(k)$, under the
connecting homomorphism
   $$\Ext^1(\Omega^1_X,F_*(\O_X)) \to
    \Ext^1(\Omega^1_X,B_1\Omega^1_X) \to
    \Ext^2(\Omega^1_X,\O_X)$$ But $F$-splitting implies that
    $\Ext^1(\Omega^1_X,B_1\Omega^1_X)$ is a direct summand of
    $\Ext^1(\Omega^1_X,F_*(\O_X))$.  This implies that the
    obstruction to lifting $X$ to a flat scheme over $W_2(k)$,
    which is an element of $\Ext^2(\Omega^1_X,\O_X)$ is zero.
\end{proof}

\begin{corollary}\label{liftability}
Let $X/k$ be a smooth, projective, Frobenius split variety. 
The $X$ admits a flat lifting to $W_2$.
\end{corollary}

\begin{remark}
I am indebted to John Millson for pointing out that Bogomolov has
proved that the deformation theory of Calabi-Yau varieties over
complex numbers is unobstructed. The above result implies that any
smooth projective, Frobenius split, Calabi-Yau variety lifts to $W_2$.
\end{remark}

\section{Unirational and Fano Threefolds}\label{unirational-and-fano}
   In this section we give application of our criterion for the
degeneration of the slope spectral sequence to smooth, projective
unirational or Fano threefolds.

	Let $X$ be a smooth projective variety of dimension $n$. We
say $X$ is Fano if $\Omega^n_X$ is the inverse of an ample line
bundle.  In characteristic zero, if $X$ is a Fano variety of dimension
$n$, then the Kodaira vanishing theorem gives $H^i(X,\O_X)=0$ for
$0<i\leq n$. This is not known in positive characteristic, without
some additional hypothesis.

Let $X$ be a smooth projective variety. We will say that $X$ is
unirational if there exists a dominant surjective, generically finite
and separable morphism $\P^n\to X$.

\begin{theorem}\label{unirat-fano}
Let $X$ be a smooth, projective, unirational or Fano 
threefold over a perfect field $k$. Then $X$ is
Hodge-Witt. 
\end{theorem}
\begin{proof}
	The proof is inspired by \cite{nygaard78}. By
 Theorem~\ref{criterion} we need to check that for $i\geq 0$,
 $H^i(X,W(\O_X))$ are of finite type when $X$ is Fano or a unirational
 variety. In both the cases we are reduced, by induction on using the
 exact sequence
\begin{equation}
0 \to W_n(\O_X) \to W_{n+1}(\O_X) \to \O_X \to 0
\end{equation}
to proving that $H^i(X,W_1(\O_X))=H^i(X,O_X)=0$ for $i>0$. 

	In \cite{nygaard78} it was shown that $H^i(X,\O_X)$ is zero
for unirational threefolds. In the case when $X$ is a Fano threefold
it was proved in \cite{shepherd-barron97} that for any smooth Fano
threefold, we have $H^i(X,\O_X)=0$ for $i>0$.
\end{proof}

\begin{remark}
	Using the universal coefficient theorem and Nygaard's result
(see \cite{nygaard78}), it is easy to see that $H^2_{\rm cris}(X/W)$
is torsion free if $X$ is a smooth, projective, unirational variety
with $\pi_1^{\rm alg}(X)=0$. When $X$ is a unirational threefold, we
can drop the assumption that $X$ is simply connected, as Nygaard has
proved that unirational threefolds are simply connected.  Hence one
deduces (via \cite{illusie79b}) that the Neron-Severi group of $X$ has
no $p$-torsion if $X$ is a unirational threefold. This is the
characteristic $p>0$ variant of the corresponding result in
characteristic zero (see \cite{serre58}).
\end{remark}

\section{Coda: Chow groups and Abel-Jacobi mappings}
   This section is a mere coaptation of the results of
 Section~\ref{degeneration} and the results of Gros and Suwa (see
 \cite{gros88a}, \cite{gros88b}, \cite{suwa88}).  We give some
 geometric applications of the degeneration of the slope spectral
 sequence of Frobenius split threefolds to Abel-Jacobi mappings.

   We need to recall a few facts about Abel-Jacobi mappings of
Bloch and Gros-Suwa (see \cite{bloch79}, \cite{gros88a}).  First
suppose that $\ell\neq p$ is a prime.  Let $\ch^i(\bX)$ be the group
of algebraic cycles of codimension $i$ on $\bX$ modulo rational
equivalence and let $\ch^i(\bX)_{\alg}$ be the group subgroup of
cycles which are algebraically equivalent to zero. Let
\begin{equation}
\lambda_i:\ch^i(\bX)_{\ell-\tors}\to H^{2i-1}(\bX,\Q_\ell/\Z_\ell(i))
\end{equation}
 be the $\ell$-adic Abel-Jacobi mapping defined by Bloch (see
\cite{bloch79}).

When $\ell=p$ we will use a variant of this map 
\begin{equation}
\lambda_i':\ch^i(\bX)_{p-\tors}\to H^{2i-1}(\bX,\Q_p/\Z_p(i))
\end{equation}
which was constructed by Gros and Suwa in \cite{gros88a}, \cite{gros88b}. 
The target of $\lambda_i'$ is defined in terms of 
logarithmic cohomology groups of $X$.

	For $\ell\neq p$ it is known that  $\lambda_2$ is
injective (see \cite{merkurev83}, \cite{colliot-thelene83}). This
injectivity assertion is also valid for the $p$-adic Abel-Jacobi
mapping $\lambda_2'$ and is proved in \cite{gros88a}.

\begin{theorem}\label{l-adic}
Let $X$ be a smooth projective variety over a perfect field $k$ of
characteristic $p>0$ and let $\ell\neq p$ be a prime.  Suppose that
for some $i\geq 1$, the cohomology group \begin{equation}
H^{2i-1}(X,W(\O_X))\tensor K\neq 0.  
\end{equation} 
Then the mapping
\begin{equation}
\lambda_i:\ch^i(\bX)_{\alg,\ell-\tors}\to H^{2i-1}(\bX,\Q_\ell/\Z_\ell(i))
\end{equation}
is not surjective.
\end{theorem}
\begin{proof}
Suwa \cite{suwa88} shows (see Corollary~3.4) that if $H^{2i-1}_{cris}(X/W)$
 has slopes outside the interval $[i-1,i]$ then the mapping $\lambda_i$ 
is not surjective. Our hypothesis shows that $X$ has a nonzero 
slope zero part  in $H^{2i-1}_{cris}(X/W)$.
\end{proof}

\begin{corollary}
Let $X$ be a smooth, projective, Frobenius split variety with
\begin{equation}
H^{2i-1}(X,W(\O_X))\tensor K\neq0.
\end{equation} 
Then 
\begin{equation}
\lambda_i:\ch^i(\bX)_{alg, \ell-\tors}\to H^{2i-1}(\bX,\Q_\ell/\Z_\ell(i))
\end{equation} 
is not surjective for any $\ell\neq p$.
\end{corollary}

\begin{corollary}\label{calabi-yau}
	If $X$ is smooth projective Calabi-Yau, Frobenius split
variety of dimension $2n-1$, for $n\geq 2$. Then
\begin{equation}
\lambda_n:\ch^n(\bX)_{\alg,\ell-\tors}\to H^{2n-1}(\bX,\Q_\ell/\Z_\ell(n))
\end{equation} 
is not surjective.
\end{corollary}

\begin{proof}
	For a Frobenius split, Calabi-Yau variety of dimension $2n-1$,
it is easy to show that $H^{2n-1}(X,W(\O_X))$ is nonzero and by
Theorem~\ref{finiteness} it is of finite type over $W$.
\end{proof}

 We say that a $p$-torsion Abelian group is of finite cotype if its
 Pontryagin dual is direct sum of finite number of copies of $\Z_p$
 and finite group.

\begin{theorem}\label{p-adic}
	Let $X$ be a smooth projective threefold over a perfect field
$k$ of characteristic $p>0$.  Assume that $H^3(X,W(\O_X))$ is of
finite type over $W$.  Then $\ch^2(\bX)_{p-\tors}$ is an Abelian
group of finite cotype.
\end{theorem}

\begin{proof}
	By the hypothesis and Theorem~\ref{criterion} we see that the
 differential $H^2(X,W\Omega_X^1)\to H^2(X,W\Omega_X^2)$ is zero. Then
 we are done by Corollary 3.7 of \cite{gros88a}.
\end{proof}

In particular we get

\begin{corollary}
 If $X$ is a smooth, projective, Frobenius split threefold over a 
perfect field $k$ of characteristic $p>0$
then $\ch^2(\bX)_{p-\tors}$ is of finite cotype.
\end{corollary}

\begin{remark}
   Let $f(x_0,x_1,x_2,x_3,x_4)=$ be a hypersurface in $\P^4$.
Assume that the coefficient of $(x_0x_1x_2x_3x_4)^{p-1}$ is non-zero
in $f^{p-1}$.  Then this hypersurface is $F$-split. If $f=0$ defines a
smooth subvariety of $\P^4$ then we can apply Theorems \ref{l-adic}
and \ref{p-adic} to the Chow groups of $f=0$.
\end{remark}

\begin{remark}
	 We can also apply this consideration to some singular
varieties. Especially to the quintic variety, given by
\begin{equation}
x_0^5+x_1^5+x_2^5+x_3^5+x_4^5=5x_0x_1x_2x_3x_4,
\end{equation}
which was first investigated by Schoen in \cite{schoen86}. This
quintic is singular with ordinary double points, but when $p\equiv
1\mod 5$, it is Frobenius split and using a criterion of
\cite{lakshmibai98} one checks easily that its desingularization is
Frobenius split. Using Theorem~\ref{criterion} we deduce that its Chow
group of codimension two cycles has $p$-torsion of finite cotype and
using the usual formula for Chow groups of blowups we deduce that the
group of codimension two cycles on Schoen variety has $p$-torsion of
finite cotype.
\end{remark}

\begin{remark}
   	These results together with the conjectures of Bloch (see
\cite{bloch-lectures}) indicate the presence of a non-trivial
filtration on Chow groups of Frobenius split varieties and the
non-degeneracy of the coniveau filtration on these varieties.  Over
the algebraic closure of a finite field one expects the Abel-Jacobi
mappings to be surjective (this has been verified by Schoen in
\cite{schoen99} under the rubric of the Tate conjecture).
\end{remark}

\begin{remark}
	Over the algebraic closure of a finite field, our results,
together with the conjectural yoga of slopes, as expounded by Bloch
(see \cite{bloch-lectures}), and Schoen's work suggests that the
Griffiths group of Frobenius split threefolds is non-trivial. We note
that we do not know how to produce explicit cycles on Frobenius split
threefolds except possibly in a small number of cases.
\end{remark}


\end{document}